\documentclass[a4paper,10pt]{article}
\usepackage[latin1]{inputenc}
\usepackage[T1]{fontenc}
\usepackage[francais]{babel}
\usepackage{amsmath}
\usepackage{amsfonts}
\usepackage{amssymb}
\usepackage{amsthm}
\usepackage{mathrsfs}
\usepackage{graphicx}
\theoremstyle{definition}
\newtheorem{definition}{\textbf{Définition}}[section]

\newtheorem{exemple}[definition]{\textbf{Exemple}}

\newtheorem{remarque}[definition]{\textbf{Remarque}}

\theoremstyle{plain}

\newtheorem{theoreme}[definition]{\textbf{Théorème}}

\newtheorem{lemme}[definition]{\textit{Lemme}}

\DeclareMathOperator{\e}{e}
\DeclareMathOperator{\li}{li}

\title{\Huge{\bf{Moments de la fonction Delta de Hooley associée à un caractère  }}}
\author{Alexandre Lartaux}

\begin{document}

\maketitle

\textbf{Abstract} 
Let $f$ be an arithmetic function, $V\geqslant 1$ a real number and
$$\Delta_V(n,f):=\sup\limits_{\substack{u \in \mathbb{R}\\ v \in [0,V]}}{\Big|\sum\limits_{\substack{d\mid n \\ \e^{u}<d\leqslant \e^{u+v}}}{f(d)}\Big|}{\rm .}$$

In \cite{B}, La Bretèche and Tenenbaum investigated weighted moments of $\Delta_1(n,f)$ where $f$ is a non principal real Dirichlet character, or the Möbius function. Answering a question of Hooley, we extend their results studying dependance  in $V$ and including the case of complex characters.  

\tableofcontents

\section{Introduction et énoncés des résultats}

Pour $n \in \mathbb{N}^*$, $f$ une fonction arithmétique, $u \in \mathbb{R}$ et $v\in [0,1]$, nous posons

\begin{align*}
\Delta(n,f,u,v):= \sum\limits_{\substack{ d\mid n \\ \e^u<d\leqslant \e^{u+v}}}{f(d)} 
\end{align*}

\begin{align*}
\Delta(n,f):=\sup\limits_ {\substack{u \in \mathbb{R} \\ v \in [0,1]}} |\Delta(n,f,u,v)| 
\end{align*}

Dans \cite{B}, La Bretèche et Tenenbaum étudient l'ordre moyen de la fonction $\Delta(n,\chi)^2$ lorsque $\chi$ est un caractère de Dirichlet réel non trivial, nous nous intéressons ici à une généralisation de leur résultat. \\ 

Pour $n \in \mathbb{N}^*$, $f$ une fonction arithmétique, et $V \in \mathbb{R}_+^*$, nous considérons
\begin{equation}
\label{eq delta}
\Delta_V(n, f):=\sup\limits_ {\substack{u \in \mathbb{R} \\ v \in [0,V]}} |\Delta(n,f,u,v)| \mbox{,} \; \; \Delta_v^*(n, f):=\sup\limits_ {u \in \mathbb{R} } |\Delta(n,f,u,v)|
\end{equation}
 
 Dans la définition de la fonction $\Delta(n,f)$, le paramètre $v$ est compris entre~$0$ et $1$. Dans celle de $\Delta_V(n,f)$, il est compris entre $0$ et $V$ où $V\geqslant 1$ et nous démontrons une estimation, uniforme en $V$, de l'ordre moyen de la fonction $\Delta_V(n,\chi)^2$. Le caractère $\chi$ est fixé et nous notons $r$ son ordre. Ainsi nos résultats sont valables pour des caractères non réels, étendant pour $V=1$ les résultats de \cite{B}.

Lorsque $r\geqslant 2$ un entier, nous notons
$$\zeta=\zeta_r:=\e\Big(\frac{1}{r}\Big) $$
\noindent
où $\e(x)=\e^{2i\pi x} {\rm .}$
\begin{definition}
Lorsque $\chi$ est un caractère de Dirichlet complexe d'ordre $r\geqslant 2$, $A>0$, $c>0$ et $\eta \in ]0, 1[$, nous disons qu'une fonction multiplicative $g$ appartient à la classe $\mathcal{M}_A(\chi,c,\eta)$ si elle est positive ou nulle et si les conditions suivantes sont vérifiées
\begin{align}
  & \forall p \in \mathcal{P} \mbox{ et } \forall \nu \geqslant 1, \; g(p^{\nu}) \leqslant A^{\nu} \label{eq def1} \\
  & \forall \varepsilon >0 \mbox{ et } \forall n \geqslant 1,\; g(n) \ll_{\varepsilon} n^{\varepsilon}\label{eq def2} \\
  &\forall k \in [0,\ldots,r-1], \; \sum\limits_{\substack{p\leqslant x \\ \chi(p)=\zeta^k}}{g(p)}=z_k \; \li(x)+ O(x\e^{-2c(\log x)^{\eta}})\label{eq def3}  
\end{align}
\noindent
pour certaines constantes $z_k=z_k(g) \geqslant 0$. \\
\end{definition}
L'article \cite{B} traite le cas où le caractère $\chi$ est réel, notre étude couvre également le cas des caractères complexes.
\begin{remarque}
Les constantes $z_k$ peuvent être calculées facilement à partir de l'estimation des sommes
$$\sum\limits_{p \leqslant x}{g(p)\chi(p)^j}$$
\noindent
pour $j \in \{0,\ldots, r-1\}$.
\end{remarque}

Lorsque $g \in \mathcal{M}_A(\chi,c,\eta)$ et $t\geqslant 1$, nous posons
$$y:=\sum\limits_{k=0}^{r-1}{z_k} \mbox{, } \; \; \; \beta_g(r,t):=\frac{1}{2^{2t-1}}\sum\limits_{k=0}^{r-1}{z_k\Big|1+\zeta^k\Big|^{2t}} {\rm .}$$
La valeur $y$ correspond donc à l'ordre moyen de $g$ sur les nombres premiers. Pour $t\geqslant 1$, nous définissons
$$\lambda(t):= \frac{1}{2\pi}\int_{-\pi}^{\pi}{|1+\e^{i\vartheta}|^{2t}{\rm d}\vartheta}=\frac{2^{2t}\Gamma(t+1/2)}{\sqrt{\pi}\Gamma(t+1)} \rm{.}$$
\noindent
Nous pouvons d'ores et déjà remarquer que $\lambda(t)\leqslant 2^{2t-1}$ pour $t\geqslant 1$ avec égalité si et seulement si $t=1$. 

\begin{exemple}
\label{ex 1}
Considérons le cas $g(n)=y^{\omega(n)}$ et $t \in \mathbb{N}^*$. Les $z_k$ valent tous $\frac{y}{r}$ et un calcul direct fournit la formule 
$$\beta(r,t):=\beta_g(r,t)=\frac{y}{2^{2t-1}}\sum\limits_{\substack{j,k=0 \\{r\mid k-j}}}^{t}{\binom{t}{k}\binom{t}{j}}{\rm .}$$
Si nous supposons de plus que $t<r$, alors
$$\beta(r,t)=\frac{y}{2^{2t-1}}\sum\limits_{k=0}^{t}{\binom{t}{k}^2}{\rm ,}$$
tandis qu'un calcul direct de l'intégrale définissant $\lambda(t)$ fournit  
$$\lambda(t)=\sum\limits_{k=0}^{t}{\binom{t}{k}^2}{\rm ,}$$
ainsi
\begin{equation}
\label{eq ex}
\beta(r,t)=\frac{y\lambda(t)}{2^{2t-1}}{\rm .}
\end{equation}
Dans ce cas, nous avons $\beta(r,t)\leqslant y$.
\end{exemple}

Lorsque $x\geqslant 16$, nous notons
$$\mathcal{L}(x)=\e^{\sqrt{\log_2(x)\log_3(x)}} {\rm .}$$
Lorsque $t \in \mathbb{R}$,
$$t^+:=\max(t,0) {\rm .}$$
Enfin, pour $x \geqslant 2$, $V\geqslant 1$, $g$ une fonction arithmétique et $t\geqslant 1$, nous notons
\begin{equation}
\label{eq theo 1}
\mathfrak{S}_{t,V}(x,\chi,g):=\sum\limits_{n \leqslant x}{g(n)\Delta_V(n,\chi)^{2t}} {\rm .}
\end{equation}

Notre objectif est de voir comment le paramètre $V$ influe sur  le comportement de l'ordre moyen de la fonction $\Delta_V(n,\chi)^2$. Nous établissons les deux théorèmes suivants, qui répondent à une question posée par le Professeur Christopher Hooley au Professeur Gérald Tenenbaum.\\

 \begin{theoreme}
\label{theo1}
Soient $A$, $c$, $\eta$ des constantes strictement positives, $t \geqslant 1$ un réel, $\chi$ un caractère de Dirichlet non principal d'ordre $r$, $g$ une fonction de $\mathcal{M}_A(\chi,c,\eta)$ et $V \geqslant 1$. Nous posons $y_0=\frac{t}{2^{2t-1}-1}$ et nous supposons que $y>0$ et $\beta_g(r,t)\leqslant y$. Alors, il existe une constante $\alpha >0$ telle que 
$$\mathfrak{S}_{t,V}(x,\chi,g)\ll x\mathcal{L}(x)^{\alpha}V^t(\log x)^{y-1+(2^{2t-1}y-y-t)^+}{\rm .}$$
\noindent
Le facteur $\log_3(x)$ qui apparaît dans la fonction $\mathcal{L}$ peut être omis si $y> y_0$
\end{theoreme}
La fonction $g$ définie dans l'exemple \ref{ex 1} vérifie les hypothèses de notre théorème dès que $t$ est un entier strictement inférieur à $r$.

\begin{remarque}
Il est également possible d'omettre le facteur $\log_3(x)$ qui apparaît dans la fonction $\mathcal{L}$ dans le cas où, $y<y_0$. Il suffit pour cela, d'adapter la méthode itérative de \cite{B}. Cependant, cette démonstration ne sera pas incluse dans l'article du fait des complications techniques qu'elle engendre.
\end{remarque}
\begin{remarque}
Pour $V=1$, $\chi$ réel, et lorsque $t$ n'est pas un entier $\geqslant 2$, nous retrouvons le résultat de \cite{B} avec la même précision. De plus, compte-tenu de la majoration triviale 
\begin{equation}
\label{triv}
\mathfrak{S}_{t,V}\leqslant ([V]+1)^{2t}\mathfrak{S}_{t,1}\rm{,}
\end{equation}
notre majoration fait apparaître un gain de l'ordre de $V^t$. Cependant, si $t$ est un entier $\geqslant 2$ et $\chi$ réel, La Bretèche et Tenenbaum obtiennent une majoration de $\mathfrak{S}_{t,1}$ avec un meilleur exposant de $\log(x)$, sous une condition plus restrictive portant sur $y$ (formule (1.6) du théorème 1.1 de \cite{B}). Mais dans ce domaine là, nous pouvons montrer que la majoration \eqref{triv} est quasiment optimale. Ainsi, nous savons qu'une majoration avec un facteur $V^t$ n'est pas vérifiée.
\end{remarque}

\begin{remarque}
Etant donné les résultats de \cite{B} et la majoration triviale \eqref{triv}, le paramètre $V^t$ n'est significatif dans la majoration \eqref{eq theo 1} que lorsque l'on a la relation supplémentaire $\sqrt{\log_2(x)\log_3(x)}=~o(\log V)$. En effet, dans le cas contraire, le facteur $V^t$ peut être supprimé car $V^t\mathcal{L}(x)^{\alpha}=O(\mathcal{L}(x)^{\alpha'})$ où $\alpha'$ est choisi suffisamment grand.
\end{remarque}

\begin{remarque}
Le paramètre $V$ influe sur l'ordre moyen de $\Delta_V(n,\chi)^2$ de manière linéaire lorsque $\chi$ est un caractère de Dirichlet d'ordre $2$. Ce résultat montre que les compensations dues aux oscillations du caractère $\chi$ dans les intervalles $]\e^u,\e^{u+V}]$ sont de nature statistique, car $\Delta_V(n,\chi)^2$ se comporte en moyenne comme $\Delta_V(n)$ où
$$\Delta_V(n):=\sup_{u \in \mathbb{R}} \sum\limits_{\substack{ d \mid n \\ \e^u<d\leqslant \e^{u+V}}}{1} {\rm .}$$
En effet, la majoration $\Delta_V(n)\leqslant V\Delta(n)$, la minoration $\Delta_V(n)\geqslant \frac{V\tau(n)}{\log n}$, valable pour $V\leqslant \log n$ ainsi que le théorème de l'article \cite{T}, montrent l'encadrement
$$xV\ll \sum\limits_{n\leqslant x}{\Delta_V(n)} \ll xV\mathcal{L}(x)^{\alpha} {\rm .}$$
\end{remarque}

\begin{theoreme}
\label{theo 2}
Soient $\chi$ un caractère de Dirichlet non principal d'ordre $r$, $y>0$, $V\geqslant 1$ et $t \geqslant 1$. Nous posons 
\begin{equation}
\label{eq U}
U:=\min(V,\log x) {\rm .}
\end{equation}
\noindent
Lorsque $x$ tend vers $+\infty$, nous avons
\begin{align}
\label{min}
\mathfrak{S}_{t,V}(x,\chi,\mu^2y^{\omega}) \gg x\Big((\log x)^{y-1}+(\log x)^{2^ty-t-1}U^{t}+(\log x)^{2^{2t}y/r-2t-1}U^{2t}\Big) {\rm .}
\end{align}
\end{theoreme}

\begin{remarque}
Grâce au premier terme du membre de droite de \eqref{min}, lorsque $t=1$ et $y\geqslant 1$, nous retrouvons le même exposant de $\log x$, ainsi que le même exposant de $V$ dans le minorant de \eqref{min} que dans le majorant de \eqref{eq theo 1}.
\end{remarque}

Soit
\begin{align}
\label{eq y_1}
y_1(t):=\frac{t}{\lambda(t)-1}{\rm .}
\end{align}
Nous énonçons deux théorèmes qui sont des versions analogues des Théorèmes \ref{theo1} et \ref{theo 2} pour la fonction $\mu$.

\begin{theoreme}
\label{theo 1 bis}
Soient $A$, $c$, $\eta$ des constantes strictement positives, $t \geqslant 1$ un réel, $g$ une fonction de $\mathcal{M}_A(c,\eta)$ et $V\geqslant 1$. Nous supposons que $y>0$. Alors, il existe une constante $\alpha >0$ telle que pour $x\geqslant 16$
$$\mathfrak{S}_{t,V}(x,\mu,g)\ll x\mathcal{L}(x)^{\alpha}V^t(\log x)^{y-1+(\lambda(t)y-y-t)^+} {\rm .}$$
\noindent
Le facteur $\log_3(x)$ qui apparaît dans la fonction $\mathcal{L}$ peut être omis si $y>y_1$.
\end{theoreme}

\begin{theoreme}
\label{theo 2 bis}
Soient $y>0$,  $V\geqslant 1$ et $t \geqslant 1$. Nous avons, lorsque $x$ tend vers $+\infty$,
\begin{align*}
  \mathfrak{S}_{t,V}(x,\mu,\mu^2y^{\omega}) 
  \gg x\Big((\log x)^{y-1}+(\log x)^{2^ty-t-1}U^{t}\Big) {\rm ,}
 \end{align*}
 \noindent
où $U$ est défini en \eqref{eq U}.
\end{theoreme}

\section{Majoration dans le cas où $v$ est fixé}

\subsection{\'Enoncés des résultats}

Lorsque $g$ et $f$ sont deux fonctions arithmétiques,  $v,t\geqslant 1$, et $x \in \mathbb{R}$, nous définissons
$$\mathfrak{S}_{t,v}^*(x,f,g):=\sum\limits_{n \leqslant x}{g(n)\Delta_v^*(n,f)^{2t}} {\rm ,}$$
\noindent
où $\Delta_v^*(n,f)$ est définie en \eqref{eq delta}. Nous nous intéressons dans cette partie à une majoration de $\mathfrak{S}_{t,v}^*(x,\chi,g)$.

\begin{theoreme}
\label{theo 3}
Sous les mêmes hypothèses que le Théorème \ref{theo1}, pour tout $v\geqslant~1$, il existe une constante $\alpha >0$ telle que 
$$\mathfrak{S}_{t,v}^*(x,\chi,g)\ll x\mathcal{L}(x)^{\alpha}v^t(\log x)^{y-1+(2^{2t-1}y-y-t)^+} {\rm .}$$
\noindent
Le facteur $\log_3(x)$ qui apparaît dans la fonction $\mathcal{L}$ peut être omis si $y> y_0$.
\end{theoreme}
Ce théorème possède également une version analogue lorsque nous remplaçons le caractère $\chi$ par la fonction $\mu$.
\begin{theoreme}
\label{theo 3 bis}
Sous les mêmes hypothèses que le Théorème \ref{theo 1 bis}, il existe une constante $\alpha>0$ telle que
$$\mathfrak{S}_{t,v}^*(x,\mu,g)\ll x\mathcal{L}(x)^{\alpha}v^t(\log x)^{y-1+(\lambda(t)y-y-t)^+} {\rm .}$$
\noindent
Le facteur $\log_3(x)$ qui apparaît dans la fonction $\mathcal{L}$ peut être omis si $y>y_1$ où $y_1$ est défini en \eqref{eq y_1}.
\end{theoreme}
\subsection{Rappels}

Dans cette section, nous énonçons les lemmes permettant de démontrer le Théorème \ref{theo 3}. Les démonstrations étant proches de celles de l'article \cite{B}, nous nous intéressons aux modifications apportées à celles-ci par l'introduction du facteur $v$ et la généralisation aux caractères complexes. \\ 

Lorsque $n \geqslant 1$, nous posons
\begin{align*}
E(n):=\min\limits_{\substack{dd'\mid n \\ d<d'}}{\log \frac{d'}{d}}\mbox{,} & &&
E^*(n):=\min{\big(1,E(n)\big)} {\rm .}
\end{align*}
\noindent
Lorsque $n \geqslant 1$, $q \geqslant 1$, et $f$ une fonction, nous posons
$$M_{q,v}^*(n,f):=\int_{\mathbb{R}}{|\Delta(n,f,u,v)|^q{\rm d} u} {\rm .}$$
\noindent
Les deux lemmes suivants étant en tout point identiques à ceux de l'article \cite{B}, nous nous contentons d'énoncer les résultats.
\begin{lemme}
\label{lemme 3.1}
Soit $v\geqslant 1$. Lorsque $n\geqslant 1$ et $q \geqslant 1$, nous avons 

\begin{equation}
\label{eq 3.1}
\Delta_v^*(n,\chi)^2 \leqslant 2^5+2^{3+2/q}E^*(n)^{-2/q}M_{2q,v}^*(n,\chi)^{1/q} {\rm .}
\end{equation}
\noindent
De plus, si $E^*(n) \geqslant \eta$ et $1 \leqslant q \leqslant b$, alors

\begin{equation}
\label{eq 3.1 bis}
M_{2b,v}^*(n,\chi)^{1/b} \leqslant 2^{8-8q/b}\eta^{\frac{2}{b}-\frac{2}{q}}M_{2q,v}^*(n,\chi)^{1/q} {\rm .}
\end{equation}
\end{lemme}

\begin{lemme}
\label{lemme 3.2}
Soit $g \in \mathcal{M}_A(c,\eta)$. Pour tout $y>0$, nous avons, uniformément en $\sigma>0$,
\begin{equation}
\label{eq 3.2 bis}
\sum\limits_{\substack{n \geqslant 1 \\ E(n) \leqslant \sigma^{3 \cdot 4^ty}}}{\frac{\mu^2(n)g(n)\tau(n)^{2t}}{n^{1+\sigma}}} \ll 1 {\rm .}
\end{equation}
\end{lemme}

Avant d'énoncer le prochain lemme, quelques notations sont nécessaires.\\ Lorsque $q,n\geqslant 1$, $w \in \mathbb{R}$ et $0\leqslant j \leqslant q$, nous définissons

$$N_{j,q,v}(n,w):=\int_{\mathbb{R}}{|\Delta(n,\chi,u,v)|^j|\Delta(n,\chi,u-w,v)|^{q-j}{\rm d} u} {\rm .}$$
\noindent
Lorsque $\chi$ est un caractère de Dirichlet, $n$ un entier strictement positif et $\vartheta$ un réel, nous posons
\begin{equation}
\label{eq tau}
\tau(n, \chi, \vartheta):= \sum\limits_{d \mid n}{\chi(d)d^{i\vartheta}} {\rm .} 
\end{equation}
\noindent
Nous définissons, pour $v\geqslant 1$,
\begin{equation}
\label{eq tau*}
 \tau_v^*(n, \chi):=\int_{0}^{\infty}{\frac{v^2}{1+\vartheta^2v^2}|\tau(n, \chi, \vartheta)|^2 {\rm d}{\vartheta}}=v\int_{0}^{\infty}{\frac{1}{1+\vartheta^2}\Big|\tau\Big(n, \chi, \frac{\vartheta}{v}\Big)\Big|^2 {\rm d}{\vartheta}}{\rm .} 
 \end{equation}

\begin{lemme}
\label{lemme 3.3}
Pour $1\leqslant j \leqslant q-1$, $g \in \mathcal{M}_A(c,\eta)$, $n\geqslant 1$, $v\geqslant 1$ et $x\geqslant 2$, nous avons 
$$\sum\limits_{p>x}{\frac{g(p)\log p}{p}N_{2j,2q,v}(n,\log p)}\ll M_{2q,v}^*(n,\chi)^{q-2/q-1}\tau_v^*(n,\chi)^{q/q-1}+R_{q,v}(n,x)$$
\noindent
où
$$R_{q,v}(n,x) \ll \e^{-c(\log x)^{\eta}}4^qM_{2q,v}(n)^{(2q-2)/(2q-1)}v^{2q/(2q-1)}\tau(n)^{2q/(2q-1)} $$
\noindent
et
\begin{equation}
\label{eq M}
M_{h,v}(n)=\int_{\mathbb{R}}{\Delta(n,u,v)^h {\rm d}u} {\rm .}
\end{equation}
\end{lemme}

\begin{proof}
Pour tout $n\geqslant 1$, nous avons
\begin{align*}
 \sum\limits_{p>x}{\frac{g(p)\log p}{p}N_{2j,2q,v}(n,\log p)} =\int_{\mathbb{R}}{\big|\Delta(n,\chi,u,v)\big|^{2j}S_{2(q-j)}{\rm d} u}  
\end{align*}
\noindent
où 
$$S_{2h}:=\sum\limits_{p>x}{\frac{g(p)\log (p)}{p}|\Delta(n,\chi,u-\log(p),v)|^{2h}} {\rm .}$$
\noindent
Un développement de la somme définissant $\Delta$ nous permet d'écrire, dans le cas où $v\geqslant \log\big(\frac{\max d_i}{\min d_i}\big)$,

$$S_{2h}=\sum\limits_{d_1,\ldots, d_{2h} \mid n}{\chi(d_1\cdots d_{h})\overline{\chi(d_{h'+1}\cdots d_{2h})}\!\!\!\!\!\!\!\!\!\!\!\!\!\!\!\!\!\!\!\!\!\!\!\!\sum\limits_{\substack{u-\log \min d_r<\log p\leqslant u+v-\log \max d_r \\ p>x}}{\!\!\!\!\!\!\!\!\!\!\!\!\!\!\!\!\!\!\!\!\!\!\!\!\frac{g(p)\log p}{p}}} {\rm .}$$
\noindent
Par sommation d'Abel, la somme intérieure vaut
$$y\int_{u-\log \min d_r}^{u+v-\log \max d_r}{1_{[0,1]}\Big(\frac{\log x}{t}\Big){\rm d}t}+O\big(\e^{-c(\log x)^{\eta}}\big)  {\rm .}$$
\noindent
En adaptant le théorème $72$ de \cite{HT}, nous obtenons, pour tout $h\geqslant 2$, l'inégalité
$$\sum\limits_{\substack{d_1,\ldots, d_{2h} \mid n \\ \log(\max d_i)-\log(\min d_i)\leqslant v}}{1}\leqslant 2^{2h}M_{2h,v}(n)$$
\noindent
où $M_{2h,v}$ est défini en \eqref{eq M}.
\noindent
Ainsi,
\begin{align*}
S_{2h} & =y\int_{\log x}^{+\infty}{|\Delta(n,\chi,u-t,v)|^{2h}{\rm d}t}+O\big(2^{2h}M_{2h,v}(n)\e^{-c(\log x)^{\eta}}\big) \\
& \leqslant y\int_{\mathbb{R}}{|\Delta(n,\chi,u-t,v)|^{2h}{\rm d}t}+O\big(2^{2h}M_{2h,v}(n)\e^{-c(\log x)^{\eta}}\big){\rm .}
\end{align*}
\noindent
Nous obtenons donc
\begin{align*}
\sum\limits_{p>x}&{\frac{g(p)\log p}{p}N_{2j,2q,v}(n,\log p)}\\
&\leqslant  yQ_j+O(2^{2(q-j)}M_{2(q-j),v}(n)M_{2j,v}(n)\e^{-c(\log x)^{\eta}})
\end{align*}
\noindent
où
$$Q_j=\int_{\mathbb{R}}{|\Delta(n,\chi,u,v)|^{2j}{\rm d}u}\int_{\mathbb{R}}{|\Delta(n,\chi,w,v)|^{2(q-j)}{\rm d}w} {\rm .}$$
\noindent
Plusieurs adaptations de l'inégalité de Hölder décrites dans le lemme 2.4 de \cite{B} impliquent
$$Q_j\ll M_{2q,v}^*(n,\chi)^{(q-2)/(q-1)}M_{2,v}^*(n,\chi)^{q/(q-1)} {\rm ,}$$
\noindent
tandis qu'une adaptation du lemme $2.3$ de  \cite{B} fournit
$$M_{2,v}^*(n,\chi)\ll \tau_v^*(n,\chi) $$
\noindent
où $\tau_v^*(n,\chi)$ est défini en \eqref{eq tau*}.
Pour traiter le reste, nous effectuons un raisonnement similaire. Le terme positif $|\Delta(n,\chi,u,v)|$ est clairement majoré par $\Delta(n,u,v)$, puis, une inégalité de Hölder avec exposants $\frac{2q-1}{2(q-j)}$ et $\frac{2q-1}{2j-1}$ montre la majoration 
\begin{align*}
& \int_{\mathbb{R}}{\Delta(n,u,v)^{2j}{\rm d}u}\\          & \leqslant \Big(\int_{\mathbb{R}}{\Delta(n,u,v){\rm d}u}\Big)^{(2(q-j))/(2q-1)}\Big(\int_{\mathbb{R}}{\Delta(n,u,v)^{2q}{\rm d}u}\Big)^{(2j-1)/(2q-1)} {\rm .}
\end{align*}
\noindent
Le même raisonnement en remplaçant $j$ par $q-j$ fournit
\begin{align*}
&M_{2j,v}(n)M_{2(q-j),v}(n) \\
&\leqslant \Big(\int_{\mathbb{R}}{\Delta(n,u,v){\rm d}u}\Big)^{2q/(2q-1)}M_{2q,v}(n)^{(2q-2)/(2q-1)} {\rm .}
\end{align*}
\noindent
Pour conclure, nous utilisons l'estimation triviale
$$\int_{\mathbb{R}}{\Delta(n,u,v){\rm d}u}\leqslant v\tau(n)$$
\noindent
valable pour $v\geqslant 1$.
\end{proof}
\begin{lemme}
\label{lemme}
Soient $A>0$, $c>0$, $\eta\in ]0,1[$, $t\geqslant 1$, $\chi$ un caractère de Dirichlet non principal d'ordre $r$ et $g \in \mathcal{M}_A(\chi,c ,\eta)$. Nous avons uniformément pour $0<|\vartheta|\leqslant 1$ et $x\geqslant 3$
\begin{align}
\label{eq 3.4}
\begin{split}
 \sum\limits_{p \leqslant x}&{\frac{g(p)}{p}|1+\chi(p)p^{i\vartheta}|^{2t}} \\
& =y\lambda(t)\log(1+|\vartheta|\log x)+\beta_g(r,t)2^{2t-1}\log\Big(\frac{\log x}{1+|\vartheta| \log x}\Big)+O(1)
\end{split}
\end{align}
\noindent
et uniformément pour $|\vartheta|> 1$
\begin{align}
\label{eq 3.4 bis}
\sum\limits_{p \leqslant x}{\frac{g(p)}{p}|1+\chi(p)p^{i\vartheta}|^{2t}}\leqslant y\lambda(t)\log_2 x+B\log_2(2+|\vartheta|)
\end{align}
\noindent
où $B$ est une constante.
\end{lemme}
\begin{proof}
Il suffit d'adapter le lemme 2.5 de \cite{B} dans le cas où le caractère $\chi$ est à valeurs complexes.
\end{proof}

\begin{lemme}
\label{lemme 4.1}
Soit $T\geqslant 1$ un réel, $A,c>0$, $\eta \in ]0,1[$ , $\chi$ un caractère de Dirichlet non principal d'ordre $r$ et $g \in \mathcal{M}_A(\chi,c,\eta)$ telle que $y>0$ et $\beta_g(r,t)\leqslant~y$. Pour tout $0<s \leqslant \frac{1}{\log T}$ et $v\geqslant 1$, nous avons
$$\sum\limits_{n\geqslant 1}{\frac{\mu^2(n)g(n)}{a_nb_n^{1+s}}\tau_v^*(n,\chi)^t}\ll v^ts^{-2^{2t-1}y}$$
\noindent
où 
$$a_n=\prod\limits_{\substack{p\mid n \\ p\leqslant T}}{p} \mbox{, } \; \;b_n=\prod\limits_{\substack{p\mid n \\ p> T}}{p} {\rm .}$$ 
\end{lemme}

\begin{proof}
\noindent
D'après l'inégalité de Hölder, nous avons
$$\tau_v^*(n,\chi)^t\ll v^{t}\int_{0}^{\infty}{\frac{1}{1+\vartheta^2}\Big|\tau\Big(n,\chi,\frac{\vartheta}{v}\Big)\Big|^{2t}{\rm d}\vartheta} {\rm .}$$

\noindent
Donc la somme que nous cherchons à estimer est 
$$\ll v^{t}\int_{0}^{\infty}{\frac{1}{1+\vartheta^2}f\Big(s,\frac{\vartheta}{v}\Big){\rm d}\vartheta}$$
\noindent
où
$$f(s, \vartheta)=\sum\limits_{n\geqslant 1}{\frac{\mu^2(n)g(n)}{a_nb_n^{1+s}}|\tau(n,\chi,\vartheta)|^{2t}} {\rm .} $$
\noindent
Pour $s>0$, cette fonction possède un développement en produit eulérien
\begin{equation}
\label{f}
f(s, \vartheta)=\prod\limits_{p\leqslant T}{\Big(1+\frac{g(p)}{p}|1+\chi(p)p^{i\vartheta}|^{2t}\Big)}\prod\limits_{p> T}{\Big(1+\frac{g(p)}{p^{1+s}}|1+\chi(p)p^{i\vartheta}|^{2t}\Big)} {\rm .}
\end{equation}
\noindent
Nous utilisons le Lemme \ref{lemme} pour estimer $f(s,\vartheta)$. Comme les équations \eqref{eq 3.4} et~\eqref{eq 3.4 bis} utilisent une borne $x$, nous devons estimer
$$\sum\limits_{p>x}{\frac{g(p)}{p^{1+s}}|1+\chi(p)p^{i\vartheta}|^{2t}} {\rm .}$$
\noindent
Or 
\begin{align*}
\sum\limits_{p>x}{\frac{g(p)}{p^{1+s}}|1+\chi(p)p^{i\vartheta}|^{2t}} & \ll 2^{2t}y \int_{x}^{\infty}{\frac{1}{t^{1+s}\log t}{\rm d}t} \\
 & \ll \int_{\log x}^{\infty}{\frac{\e^{-su}}{u}{\rm d}u} \\
 & \ll \int_{s\log x}^{\infty}{\frac{\e^{- u}}{u}{\rm d}u} {\rm .}
\end{align*}
 
\noindent
Ainsi, en prenant $x=\e^{1/s}$, cette quantité est $O(1)$. Ce choix de $x$ est compatible à la majoration de $f$ fournie par le Lemme \ref{lemme}. En effet, comme $s\leqslant \frac{1}{\log T}$, nous avons $x=\e^{1/s}\geqslant T$, ainsi, d'après \eqref{f}, nous avons 
$$f(s,\vartheta)\leqslant \prod\limits_{p\leqslant x}{\Big(1+\frac{g(p)}{p}|1+\chi(p)p^{i\vartheta}|^{2t}\Big)}\prod\limits_{p> x}{\Big(1+\frac{g(p)}{p^{1+s}}|1+\chi(p)p^{i\vartheta}|^{2t}\Big)} {\rm .}$$
\noindent
Les équations \eqref{eq 3.4} et \eqref{eq 3.4 bis} fournissent alors
$$f(s,\vartheta) \ll \left\{
 \begin{array}{ll}
\Big(1+\frac{\vartheta}{s}\Big)^{y\lambda(t)}(s+\vartheta)^{-2^{2t-1}\beta_g(r,t)} & \mbox{si }   \vartheta \leqslant 1 \\
s^{-\lambda(t)y}\log(2+\vartheta)^B & \mbox{sinon.} \\
\end{array}
\right.$$
\noindent
Ainsi
\begin{align*}
\int_{0}^{\infty}{\frac{1}{1+\vartheta^2}f\Big(s, \frac{\vartheta}{v}\Big){\rm d}\vartheta} & \ll\int_{0}^{v}{\frac{1}{1+\vartheta^2}\Big(1+\frac{\vartheta}{vs}\Big)^{\lambda(t)y}\Big(s+\frac{\vartheta}{v}\Big)^{-2^{2t-1}\beta_g(r,t)}{\rm d}\vartheta} \\
& \;\;\;\;\;+\int_{v}^{\infty}{\frac{1}{1+\vartheta^2}s^{-\lambda(t)y}\log(2+\vartheta)^B{\rm d}\vartheta} \\
& \ll s^{-2^{2t-1}\beta_g(r,t)}\int_{0}^{vs}{\frac{{\rm d}\vartheta}{1+\vartheta^2}}\\
& \;\;\;\;\;+s^{-y\lambda(t)}\int_{vs}^{v}{\Big(\frac{\vartheta}{v}\Big)^{y\lambda(t)-\beta_g(r,t)2^{2t-1}}\frac{{\rm d}\vartheta}{1+\vartheta^2}}+s^{-y\lambda(t)}\\
& \ll \min(1,vs)s^{-2^{2t-1}\beta_g(r,t)}+s^{1-2^{2t-1}\beta_g(r,t)}\log(1/s)+s^{-y\lambda(t)} {\rm .}
\end{align*}
\noindent
L'hypothèse $\beta_g(r,t)\leqslant y$ permet de conclure.
\end{proof}

Ainsi, nos hypothèses sur $\beta_g(r,t)$ sont plus restrictives que celles de $\cite{B}$, cela vient du fait que nous dilatons un intervalle portant sur les petites valeurs de $\vartheta$ d'un facteur $v$.\\ 

Avant de commencer la démonstration du Théorème \ref{theo 3}, nous énonçons les deux lemmes suivants. Le premier nous permet, de manière classique, de nous ramener au cas d'une somme portant sur les termes sans facteur carré et pondérés par un poids $1/n$. Le second correspond au lemme 70.2 de \cite{HT}.

\begin{lemme}
\label{lemme 3.5}
Soient $A>0$, $c>0$, et $\eta \in ]0,1[$. Pour tout $g \in \mathcal{M}_A(c, \eta)$, $v\geqslant 1$ et $x\geqslant 2$, nous avons 
$$\mathfrak{S}_{t,v}^*(x,\chi,g) \ll \frac{x}{\log x}\sum\limits_{n\leqslant x}{\frac{\mu^2(n)g(n)}{n}\Delta_v^*(n,\chi)^{2t}} {\rm .}$$
\end{lemme}

\begin{lemme}
\label{lemme  diff}
Soient $L(s)$, $X(s)$ deux fonctions de classe $\mathcal{C}^1$ sur $]0, \sigma_0]$, telles que 
$$-L'(s)\leqslant \phi\big(s,L(s)\big), \; \; -X'(s)\geqslant \phi\big(s,X(s)\big)$$
\noindent
où pour $s$ fixé, $\phi(s,x)$ est une fonction croissante de $x$. Si $L(\sigma_0)\leqslant X(\sigma_0)$, alors pour tout $s\leqslant \sigma_0$, $L(s)\leqslant X(s)$.
\end{lemme}

\newpage
\subsection{Démonstration du Théorème \ref{theo 3}}

La démonstration du Théorème \ref{theo 3} dans le cas $y>y_0$ est identique à celle de \cite{B}, à ceci près que nous utilisons le Lemme \ref{lemme 4.1} en lieu et place du lemme 2.7 de \cite{B}, ce qui fait apparaître un facteur $v^t$. Les Lemmes \ref{lemme 3.1}, \ref{lemme 3.2} et \ref{lemme 3.5} nous permettent de ramener notre problème à une majoration de la fonction $L_{T,q,v}$ définie par
$$L_{T,q,v}(s):=\sum\limits_{n\geqslant 1}{\frac{\mu^2(n)g(n)}{a_nb_n^{1+s}}M_{2q,v}^*(n,\chi)^{t/q}} {\rm .}$$
\noindent
Des calculs analogues à ceux de \cite{B} faisant intervenir les Lemmes \ref{lemme 3.3} et \ref{lemme 4.1} ainsi que des inégalités de Hölder fournissent
$$-L'_{T,q,v}(s)\leqslant \phi\big(s,L_{T,q,v}(s)\big)$$
\noindent
avec
$$\phi(s,x)=4^{t/q}\frac{y+as}{s}x+\frac{C_1v^{t/(q-1)}x^{(q-2)/(q-1)}}{s^{1-t/q+2^{2t-1}y/(q-1)}}+\frac{\varepsilon}{s^{\ell+1}} {\rm ,}$$
\noindent
où $\varepsilon$ et $\ell$ sont des constantes bien choisies, tandis que $C_1$ est une constante absolue. Nous cherchons un majorant sous la forme
\begin{equation}
\label{X}
X(s)=\frac{K}{s^{\gamma}}+\frac{\varepsilon}{s^{\ell}}{\rm .}
\end{equation}
La constante $\gamma$ est identique à celle de \cite{B},elle vaut $2^{2t-1}y-t+\frac{t}{q}$. En revanche, la condition sur $K$ devient $K=C_2^qv^t$ pour une constante $C_2$ suffisamment grande. Un choix pertinent de $T$, similaire à celui de \cite{B}, nous permet de négliger le second membre de \eqref{X} devant le premier.  Nous obtenons ainsi
$$L_{T,q,v}(\sigma)\ll \frac{C_2^qv^t}{\sigma^{2^{2t-1}y-t+t/q}}$$
ce qui fournit le Théorème \ref{theo 3} en spécifiant $\sigma=\frac{1}{\log x}$ et $q=\sqrt{\log(1/\sigma)}$. \\ \\

Dans le cas où $y\leqslant y_0$, la démonstration est légèrement différente. Comme dans \cite{B}, nous recherchons une fonction $X$, qui majorera notre fonction $L_{T,q,v}$ sous la forme
$$X(s):=\frac{K}{s^{\gamma}}+\frac{\varepsilon}{s^{\ell}}{\rm ,}$$
les constantes $\varepsilon$ et $\ell$ sont toujours déterminées par la fonction $\phi$, mais ici, $\gamma$ est définie par 
$$\gamma:=y+\frac{b}{q}$$
où $b$ est une constante suffisamment grande. La fonction $X$ vérifie alors les hypothèses du Lemme \ref{lemme diff} dès que 
$$K\geqslant (qC_3)^qv^t$$
où $C_3$ est une constante absolue. La nouvelle majoration de la fonction $L_{T,q,v}$ devient
$$L_{T,q,v}(\sigma)\ll \frac{(qC_3)^qv^t}{\sigma^{y+b/q}}$$
ce qui fournit le Théorème \ref{theo 3} en spécifiant $\sigma=\frac{1}{\log x}$ et $q=\sqrt{\frac{\log(1/\sigma)}{\log_2(1/\sigma)}}$. \\ \\
\begin{remarque}
Dans le cas où le caractère $\chi$ est remplacé par la fonction de Möbius, tous les lemmes cités restent valables, y compris les Lemmes \ref{lemme}
et \ref{lemme 4.1} où $r=2$, $z_0=0$ et $z_1=y$. Dans ce dernier lemme, nous pouvons par ailleurs remplacer le facteur $2^{2t-1}$ par $\lambda(t)$, ce qui permet d'adapter la méthode différentielle de manière à obtenir le Théorème \ref{theo 3 bis}.
\end{remarque}

\section{Majoration dans le cas où $v$ parcourt un intervalle de taille $V$}

\subsection{De la majoration locale à la majoration globale}
Le Théorème \ref{theo 3} nous permet de démontrer le Théorème \ref{theo1} en utilisant le lemme suivant.

\begin{lemme}
\label{lemme triv}
Soit $V\geqslant 1$, pour tout $n\geqslant 1$ et $0\leqslant \ell\leqslant 1$, nous avons 
$$\Delta_V(n,\chi)\leqslant V^{1-\ell}\Delta_{V^{\ell}}^*(n,\chi)+\Delta_{V^{\ell}}(n,\chi){\rm .}$$
\end{lemme}

\begin{proof}
L'idée de la démonstration consiste à scinder un intervalle quelconque de taille $v\leqslant V$ en union disjointe d'intervalles de taille $V^{\ell}$ et d'un intervalle de taille inférieure à $V^{\ell}$. Pour tout $u \in \mathbb{R}$ et $v \in [0,V]$, nous avons
\begin{align*}
 \Delta&(n,\chi,u,v)=\sum\limits_{\substack{d \mid n \\u<\log d\leqslant u+v}}{\chi(d)} \\
& =\sum\limits_{j=0}^{[\frac{v}{V^{\ell}}]-1}{\sum\limits_{\substack{d \mid n\\u+jV^{\ell}<\log d\leqslant u+(j+1)V^{\ell}}}{\chi(d)}}+\sum\limits_{\substack{d \mid n\\u+[\frac{v}{V^{\ell}}]V^{\ell}<\log d\leqslant u+v}}{\chi(d)} \\ 
 &=\sum\limits_{j=0}^{[\frac{v}{V^{\ell}}]-1}{\Delta(n,\chi,u+jV^{\ell},V^{\ell})}+\Delta\Big(n,\chi,u+\Big[\frac{v}{V^{\ell}}\Big]V^{\ell},v-\Big[\frac{v}{V^{\ell}}\Big]V^{\ell}\Big) {\rm ,}
 \end{align*}
 \noindent
de sorte que 
 $$\Big|\Delta(n,\chi,u,v)\Big|  \leqslant \Big[\frac{v}{V^{\ell}}\Big]\Delta_{V^{\ell}}^*(n,\chi)+\Delta_{V^{\ell}}(n,\chi){\rm .}$$
\noindent
Comme $v$ est plus petit que $V$, nous obtenons le résultat annoncé.
\end{proof}

\subsection{Démonstration du Théorème \ref{theo1}}
D'après le Lemme \ref{lemme triv}, nous avons 
$$\Delta_V(n,\chi)^{2t}\leqslant 2^{2t}\big(V^{2t(1-\ell)}\Delta_{V^{\ell}}^*(n,\chi)^{2t}+\Delta_{V^{\ell}}(n,\chi)^{2t}\big) {\rm .}$$
\noindent
Nous utilisons alors l'uniformité en $V$ de la majoration du Théorème \ref{theo 3} pour obtenir
\begin{align}
\label{eq 6}
\mathfrak{S}_{t,V}(x,g,\chi)\leqslant Ax(\log x)^{y-1}\mathcal{L}(x)^{\alpha}2^{2t}\big(V^{t(2-\ell)}+V^{2t\ell}\big){\rm .}
\end{align}
\noindent
Le facteur $V^{2t\ell}$ vient de la majoration triviale $\Delta_{V^{\ell}}(n,\chi)\leqslant V^{\ell}\Delta(n,\chi)$. En prenant $\ell=\frac{2}{3}$, nous obtenons
$$\mathfrak{S}_{t,V}(x,g,\chi)\leqslant Ax(\log x)^{y-1}\mathcal{L}(x)^{\alpha}2^{2t+1}V^{t4/3} {\rm .}$$
\noindent
où $A$ est une constante absolue.
\noindent
Par ailleurs, le facteur $V^{2t\ell}$ dans la formule \eqref{eq 6} peut désormais être remplacé par $V^{t\ell4/3}$ d'après ce qui vient d'être fait, nous déterminons alors $\ell$ tel que $2-\ell=\frac{4}{3}\ell$, soit $\ell=\frac{6}{7}$. Cela nous permet de remplacer l'exposant $\frac{4}{3}$ dans la majoration ci-dessus par $\frac{4}{3}\frac{6}{7}=\frac{8}{7}$, puis nous itérons le processus. Cependant, nous ne pouvons pas le faire indéfiniment, car à chaque itération, nous augmentons la constante. Nous constatons ainsi qu'à l'étape $n$, nous avons
\begin{equation}
\label{maj 3}
\mathfrak{S}_{t,V}(x,g,\chi)\leqslant Ax(\log x)^{y-1}\mathcal{L}(x)^{\alpha}2^{n(2t+1)}V^{tu_n}  
\end{equation}
\noindent
où $u_n$ est défini par la relation de récurrence suivante 
$$u_{n+1}=\frac{2u_n}{1+u_n}{\rm , } \; u_0=2{\rm .}$$
\noindent
En effet, nous raisonnons par récurrence sur $n$. Si l'inégalité \eqref{maj 3} est vraie au rang $n$, alors pour tout $\ell \in ]0,1[$, nous avons
$$\mathfrak{S}_{t,V}(x,g,\chi)\leqslant A2^{n(2t+1)}x(\log x)^{y-1}\mathcal{L}(x)^{\alpha}(V^{t(2-\ell)}+V^{t\ell u_n}){\rm .}  $$
\noindent
Nous posons alors $\ell$ tel que 
$$2-\ell=\ell u_n$$
\noindent
c'est-à-dire
$$\ell=\frac{2}{1+u_n}{\rm .}$$
En observant que $u_{n+1}=\ell u_n=2-\ell$, nous obtenons
$$\mathfrak{S}_{t,V}(x,g,\chi)\leqslant A2^{(n+1)(2t+1)}x(\log x)^{y-1}\mathcal{L}(x)^{\alpha}V^{t u_{n+1}} {\rm .} $$ 
\noindent
Il est clair que la suite $(u_n)_{n\geqslant 1}$ converge vers $1$ et que $u_n\leqslant 1+\frac{1}{2^n}$, ce qui fournit
$$\mathfrak{S}_{t,V}(x,g,\chi)\leqslant Ax(\log x)^{y-1}\mathcal{L}(x)^{\alpha}2^{n(2t+1)}V^{t(1+\frac{1}{2^n})} {\rm .} $$
\noindent
En choisissant $n=[\frac{\log(t\log V)}{\log 2}]+1$, nous obtenons
$$\mathfrak{S}_{t,V}(x,g,\chi)\leqslant Bx(\log x)^{y-1}\mathcal{L}(x)^{\alpha}V^{t}(\log V)^{b}  $$
\noindent
où $B$ et $b$ sont deux constantes absolues. La condition $V=O(\log x)$ montre que le facteur $(\log V)^{b}$ peut être absorbé par le terme $\mathcal{L}(x)^{\alpha}$ quitte à augmenter la valeur de $\alpha$, ce qui achève la démonstration du Théorème \ref{theo1}. \\

Cette méthode s'applique dans le cas où le caractère $\chi$ est remplacé par la fonction $\mu$. Elle permet donc de déduire le Théorème \ref{theo 1 bis} du Théorème \ref{theo 3 bis}.

\section{Minoration dans le cas où $v$ parcourt un intervalle de taille $V$}

\subsection{Une application de la transformée de Fourier}
La démonstration du Théorème \ref{theo 2} étant très semblable à celle du théorème~1.2 de \cite{B}, nous ne la rédigerons pas complètement. Nous nous contenterons de souligner les détails à modifier afin d'obtenir le Théorème \ref{theo 2}. \\ 

Pour démontrer ce résultat, nous avons besoin d'introduire l'intégrale suivante. Pour $n\geqslant 1$, $V>0$ et $f$ une fonction, nous posons 
$$M_{2,V}(n,f):=\int^{V}_{0}{\int_{\mathbb{R}}{|\Delta(n,f,u,v)|^2{\rm d}u}{\rm d}v} {\rm .}$$
\noindent
Nous pouvons alors énoncer le lemme suivant.

\begin{lemme}
\label{lemme 1}
Pour $n\geqslant 1$, $V>0$ et $\chi$ un caractère de Dirichlet, nous avons 
$$M_{2,V}(n,\chi)=\frac{V}{\pi}\int_{\mathbb{R}}{\Big(1-\frac{\sin (\vartheta V)}{\vartheta V}\Big)\frac{|\tau(n,\chi,\vartheta)|^2}{\vartheta^2}{\rm d}\vartheta}{\rm .}$$
\noindent
où $\tau(n,\chi,\vartheta)$ est défini en \eqref{eq tau}.
\end{lemme}

\begin{proof}
La démonstration de ce lemme est identique au début de la démonstration du lemme 2.3 de \cite{B}. Fixons $v \in [0,V]$ et notons 
$$\begin{array}{ccccc}
Q & : & \mathbb{R} & \to & \mathbb{C} \\
 & & u & \mapsto & \Delta(n,\chi,u,v) {\rm .}\\
\end{array}$$
\noindent
La fonction $Q$ est intégrable sur $\mathbb{R}$, nous pouvons donc calculer sa transformée de Fourier 
\begin{align*}
\widehat{Q}(\vartheta) & =\sum\limits_{d\mid n}{\chi(d)\int_{\log(d)-v}^{\log d}{\e^{-i\vartheta u} {\rm d}u}} \\
& =\sum\limits_{d\mid n}{\chi(d)d^{-i\vartheta}\e^{iv\vartheta /2}\frac{\sin(v\vartheta /2)}{\vartheta /2}} \\
& =\tau(n,\chi,-\vartheta)\e^{iv\vartheta /2}\frac{\sin(v\vartheta /2)}{\vartheta /2} {\rm .}
\end{align*}
\noindent
Comme $Q$ et $\widehat{Q}$ sont $L^2$, la formule de Plancherel assure que
$$\int_{\mathbb{R}}{Q^2(u){\rm d}u}=\frac{1}{2\pi}\int_{\mathbb{R}}{|\tau(n,\chi,\vartheta)|^2\Big(\frac{\sin(v\vartheta /2)}{\vartheta /2}\Big)^2 {\rm d}\vartheta} {\rm .}$$
\noindent
Puisque
$$\int^{V}_{0}{\sin^2(\vartheta v/2){\rm d}v}=\frac{V}{2}\Big(1-\frac{\sin (\vartheta V)}{\vartheta V}\Big) {\rm .}$$
\noindent
nous obtenons l'égalité en intégrant sur $v$ entre $0$ et $V$ l'égalité précédente.
\end{proof}

\subsection{Démonstration du Théorème \ref{theo 2}}
\begin{proof}
Plaçons-nous dans le cas où $V\leqslant \log x$. La minoration
$$\sum\limits_{n \leqslant x}{\mu^2(n)y^{\omega(n)}\Delta_V(n,\chi)^{2t}}\gg x\log (x)^{y-1}$$ 
découle de la minoration triviale $\Delta_V(n,\chi)\geqslant 1$. \\
 
La minoration 
$$\sum\limits_{n \leqslant x}{\mu^2(n)y^{\omega(n)}\Delta_V(n,\chi)^{2t}}\gg x\log (x)^{2^{2t}y/r-2t-1}V^{2t}{\rm ,}$$ 
\noindent
vient de la minoration $\Delta_V(n,\chi)\geqslant \frac{h_{\chi}(n)\tau(n)V}{1+\log n}$ où $h_{\chi}$ est la fonction multiplicative qui vaut $1$ sur les nombres premiers vérifiant $\chi(p)=1$ et $0$ sur les autres nombres premiers. \\
 
Montrons maintenant la minoration  
$$\sum\limits_{n \leqslant x}{\mu^2(n)y^{\omega(n)}\Delta_V(n,\chi)^{2t}}\gg x\log (x)^{2^ty-t-1}V^{t}{\rm .}$$
\noindent
Du Lemme \ref{lemme 1}, nous déduisons l'estimation 
$$M_{2,V}(n,\chi)\geqslant \frac{V^2}{6}\int^{1/V}_{-1/V}{\Big(1-\frac{\vartheta^2 V^2}{20}\Big)|\tau(n,\chi,\vartheta)|^2 {\rm d}\vartheta}$$
\noindent
ainsi
$$\Delta_V(n,\chi)^2\geqslant \frac{1}{7} \frac{I(n)}{\log x}$$
\noindent  
où
$$I(n)=V^2\int^{1/V}_{-1/V}{|\tau(n,\chi,\vartheta)|^2 {\rm d}\vartheta} {\rm .}$$
\noindent
Posons $w=2^{t-1}y$ et 
$$S_t(x,y)=\sum\limits_{n\leqslant x}{\mu^2(n)y^{\omega(n)}\Delta(n,\chi)^{2t}} {\rm .}$$
\noindent
D'après l'inégalité de Hölder et l'inégalité précédente, nous avons 
\begin{align*}
\sum\limits_{n \leqslant x}{\mu^2(n)w^{\omega(n)}I(n)} & \leqslant \Big(\sum\limits_{n \leqslant x}{(2^ty)^{\omega(n)}}\Big)^{1-1/t}\Big(\sum\limits_{n \leqslant x}{\mu^2(n)y^{\omega(n)}I(n)^t}\Big)^{1/t} \\
 & \ll x^{1-1/t}(\log x)^{(2^ty-1)(1-1/t)}S_t(x,y)^{1/t}\log x
\end{align*}
\noindent
d'où
$$S_t(x,y)\gg \frac{x}{(\log x)^{2^t(t-1)y+1}}\Big(V^2\int^{1/V}_{-1/V}{f(x, \vartheta) {\rm d}\vartheta}\Big)^t$$
\noindent
avec
$$ f(x, \vartheta):=\frac{1}{x}\sum\limits_{n \leqslant x}{\mu^2(n)w^{\omega(n)}|\tau(n,\chi,\vartheta)|^2} {\rm .}$$
\noindent
Si nous notons $F(s)$ la série de Dirichlet associée à $f$, un développement en produit eulérien permet d'écrire
$$F(s)=\zeta(s)^{2w}L(s+i\vartheta, \overline{\chi})^wL(s-i\vartheta, \chi)^wH(s,\vartheta)$$
\noindent
avec $H(s, \vartheta)$ uniformément convergente et bornée dans le demi plan $\mathfrak{Re}(s) \geqslant \frac{3}{4}$ et vérifiant $|H(1, \vartheta)|\gg 1$ lorsque $|\vartheta| \leqslant 1$. Ainsi, d'après le théorème de Selberg et Delange \cite{TAN}, nous avons, dès que $\frac{C}{\log x} \leqslant \vartheta\leqslant \frac{1}{V}$ (pour une constante $C$ suffisamment grande)
$$f(x, \vartheta) \gg (\log x)^{2w-1} {\rm .}$$
\noindent
Nous obtenons ainsi
$$S_t(x,y)\gg x(\log x)^{2^ty-t-1}V^t {\rm .}$$
\end{proof}

\begin{remarque}
Dans le cas où le caractère $\chi$ est remplacé par la fonction~$\mu$, le seul point de la démonstration qui diffère est que la fonction $L(s,\chi)$ est remplacée par l'inverse de la fonction $\zeta$ de Riemann. La méthode de Selberg-Delange s'applique également à ce cas, ce qui nous permet de démontrer le Théorème \ref{theo 2 bis}. 
\end{remarque}

\section*{Remerciements}

Je tiens ici à remercier mon directeur de thèse, Régis de la Bretèche, et Gérald Tenenbaum, pour  leurs conseils et leurs remarques pertinentes qui m'ont été très précieux lors de l'élaboration de cet article.

\nocite{*}
\bibliographystyle{amsplain}


\bibliography{biblio4}
\addcontentsline{toc}{section}{Références}

\vglue0.3cm
\hglue0.02\linewidth\begin{minipage}{0.9\linewidth}
\begin{center}
{Universit\'e de Paris, Sorbonne Université}\\
CNRS,  \\
 Institut de~Math\'ematiques~de Jussieu- Paris Rive Gauche,\\
F-75013 Paris, France \\
 E-mail : \parbox[t]{0.45\linewidth}{\texttt{alexandre.lartaux@imj-prg.fr}}

\end{center}
\end{minipage}
\end{document}